\documentclass[12pt]{amsart}
\usepackage{amsmath}
\usepackage{amssymb, amsfonts, amsthm}
\usepackage{amssymb}
\usepackage{ifthen}
\usepackage{graphicx}
\usepackage{float}
\usepackage{caption}
\usepackage{subcaption}
\usepackage{cite}
\usepackage{amsfonts}
\usepackage{amscd}
\usepackage{amsxtra}
\usepackage{color}
\numberwithin{equation}{section}

\newtheorem{theorem}{Theorem}[section]

\newtheorem{corollary}[theorem]{Corollary}

\theoremstyle{definition}

\newtheorem{remark}[equation]{Remark}

\begin{document}
\title[Bohr inequalities for certain integral operators]
{Bohr inequalities for certain integral operators}

\author[Shankey Kumar]{Shankey Kumar}
\address{Discipline of Mathematics \\ Indian Institute of Technology Indore\\
Simrol, Khandwa Road\\
Indore 453 552, India}
\email{shankeygarg93@gmail.com}

\author[S. K. Sahoo]{Swadesh Kumar Sahoo}
\address{Discipline of Mathematics \\ Indian Institute of Technology Indore\\
Simrol, Khandwa Road\\
Indore 453 552, India}
\email{swadesh.sahoo@iiti.ac.in}


\begin{abstract}
In this article, we determine sharp Bohr-type radii for certain complex integral operators 
defined on a set of bounded analytic functions in the unit disk.
\end{abstract}

\maketitle

\section{\bf Introduction }
We denote by $\mathbb{D}:=\{z\in \mathbb{C}:|z|<1 \}$ the unit disk in the complex plane.
Let $\mathcal{H}$ be the class of all analytic functions defined on $\mathbb{D}$. 
Setting $\mathcal{B}=\{f\in\mathcal{H}: |f(z)|\leq 1\}$.  
Let us first highlight a remarkable 
result of Bohr \cite{Bohr14} that opens up a new type of research problems in
geometric function theory, which states that
``{\em If $f(z)=\sum_{n=0}^{\infty} a_n z^n \in \mathcal{B}$, then
$$
\sum_{n=0}^{\infty}|a_n|r^n\leq 1,
$$
for $r\leq 1/3$ and the constant $1/3$ cannot be improved.}''
The quantity $1/3$ is known as the {\em Bohr radius} for the class $\mathcal{B}$.
Moreover, for functions in 
$\mathcal{B}_0 := \{f\in \mathcal{B} \ | \ f(0)=0 \}$,
Bombieri \cite{Bombieri62} found the Bohr radius, 
which is $1/\sqrt{2}$ (for more generalization of this result see \cite{PW2020}).  
These are not the only classes of the analytic functions where the Bohr radii are studied
but also for many other classes of functions and for some integral operators. Some of those
are highlighted below. In fact an interesting application of Bohr radius problem for 
the class $\mathcal{B}$ can be found in \cite{PVW201911}.

In \cite{Ali17}, Ali et al. brought into the notice of the Bohr radius problem for the 
odd analytic functions, which is settled by Kayumov and Ponnusamy in \cite{Kayumov17}. 
Also, Kayumov and Ponnusamy \cite{Kayponn18} generalized the problem of the Bohr radius 
for the odd analytic functions. Bhowmik and Das \cite{BD18} studied the Bohr radius for 
families of certain analytic univalent (one-to-one) functions. In \cite{BDK04}, 
the Bohr phenomenon is discussed for the functions in Hardy spaces. The study of the Bohr 
radius of the Bloch functions discussed in \cite{KSS2017}. In \cite{AlKayPON-19,PVW2019}, 
authors  studied the Bohr phenomenon for
a quasi-subordination family of functions. 
Recently, Bhowmik and Das \cite{Bhowmik19} studied the Bohr radius 
for derivatives of analytic functions.
To find more achievements in this context, 
one may see the papers \cite{Abu,Abu2,Abu4,AAL20,BD19,Kaypon18,Kayumov18,Liu18,LLP2020,LPW2020,LP2019} 
and the references there in. Also, the survey article 
\cite{Abu-M16} and the references cited in it are useful in this direction.

A natural question arises {\em``can we find Bohr radius for certain 
complex integral operators defined either on the class $\mathcal{B}$ or $\mathcal{B}_0$?"}. 
This idea has been initiated first for the classical Ces\'{a}ro operator in 
\cite{Kayponn19}. As our results of this paper 
are motivated by \cite{Kayponn19}, here first we recall the definition of 
the Ces\'{a}ro operator
followed by 
statement of the result on absolute sum of the series representation of the operator.
The Ces\'{a}ro operator is studied in \cite{Hardy32} (see, for more information, 
\cite{Siskakis87} and \cite{Siskakis90}) which is defined as
\begin{equation}\label{5eq1.1}
T[f](z):=\int_{0}^{1} \frac{f(tz)}{1-tz} dt =\sum_{n=0}^{\infty}\left(\frac{1}{n+1}\sum_{k=0}^{n} 
a_k\right) z^n,  
\end{equation}
where $f(z)=\sum_{n=0}^{\infty}a_n z^n$ is analytic in $\mathbb{D}$. 
Also, a generalized form of the Ces\'{a}ro operator is studied in \cite{Agrawal05}.

As noted in \cite{Kayponn19},
$$
|T[f](z)|\leq \frac{1}{r} \log \frac{1}{1-r}
$$
for each $|z|=r<1$. 
On the other hand, from \eqref{5eq1.1}, we also have the obvious estimate
$$
|T[f](z)|\leq \sum_{n=0}^{\infty}  \bigg(\frac{1}{n+1} \sum_{k=0}^{n} 
|a_k| \bigg)|z|^n,
$$
the absolute sum of the series \eqref{5eq1.1}.
However, if $|z|=r<1$, Kayumov et al. \cite{Kayponn19} obtained the sharp radius $r$ for which 
this absolute sum has the same upper bound $(1/r)\log(1/(1-r))$. 
This was important to study, as in general, a convergent series need not be absolutely convergent.
Indeed, they established

\medskip
\noindent
{\bf Theorem A.} {\em If $f(z)=\sum_{n=0}^{\infty}a_n z^n\in \mathcal{B}$, then
$$
\sum_{n=0}^{\infty}  \bigg(\frac{1}{n+1} \sum_{k=0}^{n} 
|a_k| \bigg)r^n\leq \frac{1}{r} \log \frac{1}{1-r}
$$
for $r\leq R= 0.5335\ldots$. Here the number $R$ is the positive root of the equation
$$
2x- 3(1-x) \log \frac{1}{1-x}=0
$$
that cannot be improved.}

\medskip
Motivated by Theorem A, in this paper, we study the Bohr radius problem for the {\em $\beta$-Ces\'aro 
operator} ($\beta>0$) defined by
$$
T_\beta[f](z):=\sum_{n=0}^{\infty} \bigg(\frac{1}{n+1}\sum_{k=0}^{n} \frac{\Gamma{(n-k+\beta)}}
{\Gamma{(n-k+1)}\Gamma(\beta)}a_{k}\bigg) z^n
=\int_{0}^{1} \frac{f(tz)}{(1-tz)^\beta} dt, \ z\in \mathbb{D},
$$
and for the Bernardi operator defined as
$$
L_\gamma [g](z):= \sum_{n=m}^{\infty} \frac{a_n}{n+\gamma} z^n
= \int_{0}^{1} g(zt) t^{\gamma-1} dt,
$$
for $g(z)=\sum_{n=m}^{\infty} a_n z^n$ and $\gamma > -m$, 
here $m\geq 0$ is an integer.  With the help of the Bernardi operator 
we also obtain the Bohr radii for some known operators.
Detailed discussion on these problems are described in the next section.
\section{\bf Main results}

Note that the $\beta$-Ces\'aro operator $T_\beta$ ($\beta>0$) is a natural 
generalization of the Ces\'{a}ro operator $T$ defined by \eqref{5eq1.1} and 
indeed, we have $T_1=T$.
For $f\in \mathcal{B}$ and $\beta>0$, an elementary estimation of the integral in absolute value 
gives us the sharp inequality
$$ 
|T_\beta[g](z)|\leq 
\begin{cases}\cfrac{1}{r}\bigg[ \cfrac{1-(1-r)^{1-\beta}}{1-\beta}\bigg], & \text{if } \beta\neq 1,\\[5mm] 
\cfrac{1}{r}\log \cfrac{1}{1-r}, & \text{ if } \beta =1, 
\end{cases}
$$
for each $|z|=r<1$. In this line, similar to Theorem A, our first main result is the following.

\begin{theorem}\label{5theorem1.1}
For $f(z)=\sum_{n=0}^{\infty}a_n z^n \in\mathcal{B}$ and $0<\beta \neq 1$, we have
$$
\sum_{n=0}^{\infty} \bigg(\frac{1}{n+1}\sum_{k=0}^{n} \frac{\Gamma{(n-k+\beta)}}
{\Gamma{(n-k+1)}\Gamma(\beta)}|a_{k}|\bigg) r^n
\leq \frac{1}{r}\bigg[\frac{1-(1-r)^{1-\beta}}{1-\beta}\bigg],
$$
for $r\leq R(\beta)$, where $R(\beta)$ is the positive root of the equation
$$
\frac{3[1-(1-x)^{1-\beta}]}{1-\beta} - \frac{2[(1-x)^{-\beta}-1]}{\beta} =0.
$$
The radius $R(\beta)$ cannot be improved.
\end{theorem}

Here, it is easy to observe that if we take the limit $\beta\to 1$ 
in Theorem \ref{5theorem1.1} then we can obtain Theorem A.

\begin{remark}
Another form of the $\beta$-Ces\'{a}ro operator of a normalized analytic function 
$g(z)=\sum_{n=1}^{\infty}b_n z^n$ in $\mathbb{D}$ has been studied in the literature (see \cite{Skumar20}):
$$
C_\beta[g](z)=\int_0^1 \frac{g(tz)}{t(1-tz)^\beta} dt=\sum_{n=0}^{\infty} \bigg(\frac{1}{n+1} \sum_{k=0}^{n} 
\frac{\Gamma{(n-k+\beta)}}{\Gamma{(n-k+1)}\Gamma(\beta)}b_{k+1}\bigg)z^{n+1},
\ z\in \mathbb{D},
$$
for $\beta> 0$.
This version of the $\beta$-Ces\'{a}ro operator was initially considered to study its boundedness, 
compactness, and spectral properties, and more recently its univalency properties were 
investigated in \cite{Kumar20}.
To study its Bohr radius problem, it is necessary for us to assume that 
$g(z)=\sum_{n=1}^{\infty}b_n z^n \in\mathcal{B}_0$. 
An easy calculation  
gives us the sharp inequality, for $g\in \mathcal{B}_0$ and $\beta>0$,
$$ 
|C_\beta[g](z)|\leq 
\begin{cases} \cfrac{1-(1-r)^{1-\beta}}{1-\beta}, & \text{if } \beta\neq 1,\\ 
\log \cfrac{1}{1-r}, & \text{ if } \beta =1, 
\end{cases}
$$
for each $|z|=r<1$.
It is well-known by the Schwarz lemma that if $g(z)=\sum_{n=1}^{\infty}b_n z^n \in\mathcal{B}_0$ 
then we can write $g(z)=zh(z)$ for $h(z)=\sum_{n=0}^{\infty}b_{n+1} z^n \in\mathcal{B}$. 
So, we have
$$
C_\beta[g](z)=\sum_{n=0}^{\infty} \bigg(\frac{1}{n+1} \sum_{k=0}^{n} 
\frac{\Gamma{(n-k+\beta)}}{\Gamma{(n-k+1)}\Gamma(\beta)}b_{k+1}\bigg)z^{n+1}
=zT_\beta[h](z).
$$
Now, by using Theorem \ref{5theorem1.1} we obtain 
$$
\sum_{n=0}^{\infty} \bigg(\frac{1}{n+1} \sum_{k=0}^{n} 
\frac{\Gamma{(n-k+\beta)}}{\Gamma{(n-k+1)}\Gamma(\beta)}|b_{k+1}|\bigg)r^{n+1}\leq\frac{1-(1-r)^{1-\beta}}{1-\beta}, \ 0<\beta\neq 1,
$$
for $r\leq R(\beta)$. Here $R(\beta)$ is the positive root of the equation
$$
\frac{3[1-(1-x)^{1-\beta}]}{1-\beta} - \frac{2[(1-x)^{-\beta}-1]}{\beta} =0
$$
that cannot be improved. 
Recall that the operator $C_1$ has been considered in \cite{Hartmann74, Ponnusamy19,Kumar20,Skumar20}
for various aspects. Moreover, in the limit $\beta\to 1$,
we can indeed obtain the Bohr radius problem:
{\em If $g(z)=\sum_{n=1}^{\infty}b_n z^n \in\mathcal{B}_0$ then 
$$
\sum_{n=0}^{\infty} \bigg(\frac{1}{n+1}\sum_{k=0}^{n} |b_{k+1}| \bigg)r^{n+1}\leq \log \frac{1}{1-r}
$$
for $r\leq R= 0.5335\ldots$. The number $R$ is the positive root of the equation
$$
2x- 3(1-x) \log \frac{1}{1-x}=0
$$
that cannot be improved.}
This remark observes that the Bohr radii for the operators $T_\beta$ and $C_\beta$ are almost same,
but up to an extra factor $1/r$.
\hfill{$\Box$}
\end{remark}

Similar to the Bohr-type radius problem of the operator $T_\beta$, $\beta>0$, 
we also study the Bohr radius of the absolute series of the Bernardi operator 
\cite[P. 11]{MM-Book} (see also \cite{PS08}) defined by
$$
L_\gamma [f](z):= \sum_{n=m}^{\infty} \frac{a_n}{n+\gamma} z^n
= \int_{0}^{1} f(zt) t^{\gamma-1} dt,
$$
for $f(z)=\sum_{n=m}^{\infty} a_n z^n$ and $\gamma > -m$, 
here $m\geq 0$ is an integer. The function $L_\gamma[f]$ is analytic 
in $\mathbb{D}$ and several properties of $L_\gamma[f]$ when $m=1$ 
(with a normalization) are well-known (see, for instance \cite{MM-Book,Parvatham08,PS08}).

It is easy to calculate the following sharp bound
$$
|L_\gamma [f](z)|\leq \frac{1}{m+\gamma}r^m, \ |z|=r<1
$$
for $f(z)=\sum_{n= m}^{\infty} a_n z^n$. Corresponding to the above inequality, we obtain the following result.

\begin{theorem}\label{5theorem1.2}
Let $\gamma>-m$. If $f(z)=\sum_{n=m}^{\infty} a_n z^n \in \mathcal{B}$, then
$$
 \sum_{n=m}^{\infty} \frac{|a_n|}{n+\gamma}r^{n}\leq \frac{1}{m+\gamma}r^m
$$
for $r\leq R(\gamma)$. Here, $R(\gamma)$ is the positive root of the equation
$$
\frac{x^m}{m+\gamma}-2
\sum_{n=m+1}^{\infty}\frac{x^n}{n+\gamma}=0
$$
that cannot be improved.
\end{theorem}

Letting $\gamma=1$ in the Bernardi operator $L_\gamma$,
we obtain the well-known Libera operator \cite{MM-Book,PS08} defined as
$$
L[f](z):=\int_{0}^{1} f(zt)\,dt=\sum_{n=0}^{\infty} \frac{a_n}{n+1} z^n.
$$
The multiplication of $z$ in the Libera operator $L$ gives the integral
$$
I[f](z):=
\sum_{n=0}^{\infty} \frac{a_n}{n+1} z^{n+1} = \int_{0}^{z} f(w) dw ,\ |z|<1. 
$$
It is easy to check that
$$
|L[f](z)|\leq 1 \ \text{and}\  |I[f](z)|\leq r, \ |z|=r.
$$ 
As a special case of Theorem \ref{5theorem1.1} ($\gamma=1$ and $m=0$), 
we get the Bohr radius for the Libera operator as well as for the operator $I$ as follows.

\begin{corollary}\label{5corollary1.1}
If $f(z)=\sum_{n=0}^{\infty} a_n z^n \in \mathcal{B}$, then
$$
\sum_{n=0}^{\infty} \frac{|a_n|}{n+1} r^n\leq 1,
$$
for $r\leq R$ with $R=0.5828\ldots$, the positive root of the equation 
$3x+2\log(1-x)=0$. Here, $R$ is the best possible.
\end{corollary}

Also, the Alexander operator \cite{Duren83, MM-Book, Kumar20, Skumar20}
$$
J[g](z):=\int_{0}^{1} \frac{g(zt)}{t}dt=\sum_{n=1}^{\infty}\frac{ b_n}{n} z^n,
$$
for $g(z)=\sum_{n=1}^{\infty} b_n z^n$, extensively studied in the univalent function theory. 
We have sharp bound
$$
|J[g](z)|\leq r
$$
for each $|z|=r<1$, since $|g(zt)/t|\le 1$ here.
Then from the observation of the Schwarz lemma, for every $g\in\mathcal{B}_0$ 
we can obtain an element $h\in \mathcal{B}$ such that $g(z)=zh(z)$. 
So, we have the following result as a consequence of Corollary \ref{5corollary1.1}.

\begin{corollary}
If $g(z)=\sum_{n=1}^{\infty} b_n z^n \in \mathcal{B}_0$, then
$$
\sum_{n=1}^{\infty} \frac{|b_n|}{n} r^n\leq r,
$$
for $r\leq R$. Here, $R=0.5828\cdots$ is the positive root of the equation 
$3x+2\log(1-x)=0$ that cannot be improved.
\end{corollary}
In the next section, we discuss the proofs of Theorems \ref{5theorem1.1} and \ref{5theorem1.2}.

%%%%%%%%%%%%%%%%%%%%%%%%%%%%%%%%%%%%%%%%%%
%%%%%%%%%%%%%%%%%%%%%%%%%%%%%%%%%%%%%%%%%

\section{\bf Proofs of the main results}

%%%%%%%%%%%%%%%%%%%%%%%%%%%%%%%%%%%%%%
%%%%%%%%%%%%%%%%%%%%%%%%%%%%%%%%%%%%%%%%%

\subsection {\bf Proof of Theorem \ref{5theorem1.1}}
First we define
\begin{equation}\label{5eq2.1}
T_\beta^f (r):= \sum_{n=0}^{\infty}\bigg(\frac{1}{n+1} \sum_{k=0}^{n} \frac{\Gamma{(n-k+\beta)}}{\Gamma{(n-k+1)}
\Gamma(\beta)}|a_{k}|\bigg)r^n,
\end{equation}
where $f(z)=\sum_{n=0}^{\infty} a_n z^n \in \mathcal{B}$, $0<\beta\neq 1$ and $r=|z|<1$. Setting $|a_0|:= a$ and let $a<1$. 
By Wiener's estimate we know that $|a_n|\leq 1-a^2$ for $n\geq 1$. This yields
$$
T_\beta^f (r) \leq a\sum_{n=0}^{\infty}\bigg(\frac{1}{n+1} \frac{\Gamma{(n+\beta)}}{\Gamma{(n+1)}\Gamma(\beta)}
\bigg)r^n+ (1-a^2)\sum_{n=1}^{\infty} \bigg(\frac{1}{n+1}\sum_{k=1}^{n} \frac{\Gamma{(n-k+\beta)}}{\Gamma{(n-k+1)}
\Gamma(\beta)}\bigg)r^n.
$$
The above inequality is equivalent to
\begin{align*}
T_\beta^f (r)
& \leq \frac{a}{r}\int_0^r \frac{1}{(1-t)^\beta}\,dt + \frac{(1-a^2)}{r}\int_0^r \frac{t}{(1-t)^{\beta+1}}\,dt\\
& = \frac{(a^2+a-1)}{r}\int_0^r \frac{1}{(1-t)^\beta}\,dt + \frac{(1-a^2)}{r}\int_0^r \frac{1}{(1-t)^{\beta+1}}\,dt.
\end{align*}
It follows that
$$
T_\beta^f(r) \leq \frac{1}{r}\bigg[ \frac{(a^2+a-1)[1-(1-r)^{1-\beta}]}{1-\beta} 
+\frac{ (1-a^2)[(1-r)^{-\beta} -1]}{\beta}\bigg]:=\phi(a).
$$
Differentiation of the function $\phi$ with respect to $a$ gives us
$$
\phi'(a)=\frac{1}{r}\bigg[\frac{(2a+1)[1-(1-r)^{1-\beta}]}{1-\beta} -\frac{ 2a[(1-r)^{-\beta} -1]}{\beta}\bigg]
$$
and so
$$
\phi''(a)=\frac{1}{r}\bigg[\frac{2[1-(1-r)^{1-\beta}]}{1-\beta} -\frac{ 2[(1-r)^{-\beta} -1]}{\beta}\bigg].
$$
It is easy to see that $\phi''(a)\leq 0$ for every $a\in[0,1)$ and $r\in[0,1)$. 
This provides that $\phi'(a)\geq \phi'(1)$. Here
$$
\phi'(1)= \frac{1}{r}\bigg[\frac{3[1-(1-r)^{1-\beta}]}{1-\beta} -\frac{ 2[(1-r)^{-\beta} -1]}{\beta}\bigg]\geq 0
$$
holds for $r\leq R(\beta)$, where
$R(\beta)$ is the positive root of the equation
$$
\frac{3[1-(1-x)^{1-\beta}]}{1-\beta} - \frac{2[(1-x)^{-\beta}-1]}{\beta} =0.
$$
Then $\phi(a)$ is an increasing function of $a$, for $r\leq R(\beta)$. 
It implies that
$$
\phi(a)\leq\phi(1)= \frac{1}{r}\bigg[\frac{1-(1-r)^{1-\beta}}{1-\beta}\bigg],
$$
for $r\leq R(\beta)$. It is easy to observe that $R(\beta)< 1$. This completes the first part of the theorem.

To conclude the final part, we consider the function
$$
\phi_a(z)= \frac{z-a}{1-az}=-a+(1-a^2)\sum_{n=1}^{\infty} a^{n-1} z^{n},
$$
where $z\in\mathbb{D}$ and $a\in [0,1)$. By using \eqref{5eq2.1}, 
we obtain the sum 
\begin{align*}
T_\beta^{\phi_{a}}(r)
& = \frac{a}{r}\bigg[\frac{1-(1-r)^{1-\beta}}{1-\beta}\bigg]+ (1-a^2)\sum_{n=1}^{\infty} \bigg(\frac{a^{n-1}}{n+1}\sum_{k=1}^n\frac{\Gamma(n-k+\beta)}{\Gamma(n-k+1)\Gamma(\beta)}\bigg)
r^n\\
& = \frac{a}{r}\bigg[\frac{1-(1-r)^{1-\beta}}{1-\beta}\bigg] + \frac{(1-a^2)}{r}\int_{0}^{r} 
\frac{t}{(1-at)(1-t)^\beta} dt.
\end{align*}
We can rewrite the last expression as
\begin{equation}\label{5eq2.2}
T_\beta^{\phi_a}(r)= \frac{1}{r}\bigg[\frac{1-(1-r)^{1-\beta}}{1-\beta}\bigg] -\frac{(1-a)}{r}\Bigg[\frac{3[1-(1-r)^{1-\beta}]}{1-\beta} 
-\frac{ 2[(1-r)^{-\beta} -1]}{\beta}\Bigg] + N_a(r),
\end{equation}
where
$$
N_a(r)= \frac{2(1-a)}{r}\Bigg[\frac{[1-(1-r)^{1-\beta}]}{1-\beta} -\frac{ [(1-r)^{-\beta} -1]}{\beta}\Bigg] 
+\frac{(1-a^2)}{r}\int_{0}^{r} \frac{t}{(1-at)(1-t)^\beta} dt.
$$
Expressing $N_a(r)$ into its summation form, we have
\begin{align*}
N_a(r)
=& \sum_{n=0}^{\infty} \frac{1}{n+1}\Bigg(-\frac{(1-a)^2}{a} 
\frac{\Gamma{(n+\beta)}}{\Gamma{(n+1)}\Gamma(\beta)}
-2(1-a)\frac{\Gamma{(n+\beta+1)}}{\Gamma{(n+1)}\Gamma{(\beta+1)}}\\
& \hspace*{3cm}+\frac{(1-a^2)}{a} \sum_{m=0}^{n}\frac{\Gamma{(n-m+\beta)}}{\Gamma{(n-m+1)}\Gamma(\beta)}a^m \Bigg) 
r^n.
\end{align*}
By using the identity
$$
\sum_{m=0}^{n}\frac{\Gamma{(n-m+\beta)}}{\Gamma{(n-m+1)}\Gamma(\beta)}
=\frac{\Gamma{(n+\beta+1)}}{\Gamma{(n+1)}\Gamma{(\beta+1)}},
$$
we can get that
$N_a(r)=O((1-a)^2)$, as $a$ tends to $1$.
Further, a simple computation shows that for $r>R(\beta)$ the quantity
$$
\frac{3[1-(1-r)^{1-\beta}]}{1-\beta} -\frac{ 2[(1-r)^{-\beta} -1]}{\beta}<0.
$$
After using these observations in \eqref{5eq2.2} we conclude that $R(\beta)$ 
cannot be improved. This completes the proof.
\hfill{$\Box$}

%%%%%%%%%%%%%%%%%%%%%%%%%%%%%%%%
%%%%%%%%%%%%%%%%%%%%%%%%%%%%%%%%
%%%%%%%%%%%%%%%%%%%%%%%%%%%%%%%%
%%%%%%%%%%%%%%%%%%%%%%%%%%%%

\subsection {\bf Proof of Theorem \ref{5theorem1.2}}
Given that $f(z)=\sum_{n= m}^{\infty}a_n z^n\in\mathcal{B}$. We set the notation
\begin{equation}\label{5eq2.3}
L_f (r):= \sum_{n= m}^{\infty}\frac{|a_{n}|}{n+\gamma}r^n.
\end{equation}
The Schwarz lemma gives $f(z)=z^m h(z)$, where $h(z)=\sum_{n=m}^{\infty}a_{n} z^{n-m}$. 
Denoting by $a:=|a_m|<1$ and using the Wiener estimate $|a_n|\leq (1-a^2)$ for $n\geq m+1$
in \eqref{5eq2.3}, we obtain the following inequality
$$
L_f (r) \leq \frac{a}{m+\gamma}r^m+ (1-a^2)\sum_{n= m+1}^{\infty}\frac{1}{n+\gamma}r^n:= \psi(a).
$$
It is easy to see that
$$
\psi''(a)= -2\sum_{n=m+1}^{\infty}\frac{1}{n+\gamma}r^n\leq 0.
$$
Thus,
$$
\psi'(a)\geq \psi'(1) = \frac{1}{m+\gamma}r^m-2\sum_{n= m+1}^{\infty}\frac{1}{n+\gamma}r^n \geq 0,
$$
for $r\leq R(\gamma)$, where
$R(\gamma)$ is the positive root of the equation
$$
\frac{1}{m+\gamma}r^m-2\sum_{n=m+1}^{\infty}\frac{1}{n+\gamma}r^n=0.
$$
Hence, $\psi(a)$ is an increasing function of $a$ for $r\leq R(\gamma)$.
This gives that
$$
\sum_{n=m}^{\infty}\frac{|a_{n}|}{n+\gamma}r^n
\leq \frac{1}{m+\gamma}r^m, \ \text{ for } r\leq R(\gamma).
$$
Also, a simple observation gives $R(\gamma)<1$.

To prove $R(\gamma)$ is the best possible bound, we consider the function
$$
\psi_a(z)=z^m \frac{z-a}{1-az}=-az^m+(1-a^2)\sum_{n=1}^{\infty} a^{n-1} z^{n+m},
$$
where $z\in\mathbb{D}$ and $a\in [0,1)$. We obtain the following equality
$$
L_{\psi_a} (r)= \frac{a}{m+\gamma}r^m+ (1-a^2) 
\sum_{n=m+1}^{\infty} \frac{a^{n-1}}{n+\gamma}r^n
$$
with the help of \eqref{5eq2.3}, which is equivalent to
\begin{equation}\label{5eq2.4}
L_{\psi_a} (r) =\frac{1}{m+\gamma}r^m- (1-a)\Bigg( \frac{1}{m+\gamma}r^m
-2
\sum_{n=m+1}^{\infty}\frac{1}{n+\gamma}r^n\Bigg) + M_a(r),
\end{equation}
where
$$
M_a(r)= 2(a-1)
\sum_{n=m+1}^{\infty}\frac{1}{n+\gamma}r^n
 + (1-a^2) 
\sum_{n=m+1}^{\infty}  \frac{a^{n-1}}{n+\gamma}r^n.
$$
Letting $a\to 1$, we obtain
$$
M_a(r)=\sum_{n=m+1}^{\infty}\frac{2(a-1)+(1-a^2) a^{n-1}}{n+\gamma} 
r^{n}= O((1-a)^2).
$$
Further, the quantity
$$
\frac{1}{m+\gamma}r^m
-2
\sum_{n=m+1}^{\infty}\frac{1}{n+\gamma}r^n<0
$$
whenever $r> R(\gamma)$. These facts in \eqref{5eq2.4} gives that $R(\gamma)$ 
cannot be improved and the proof is complete.

\bigskip
\noindent
{\bf Acknowledgement.} 
The authors thank Professor S. Ponnusamy for bringing some of the Bohr radius papers including \cite{Kayponn19} 
to their attention and useful discussion on this topic.
The work of the first author is supported by CSIR, New Delhi (Grant No: 09/1022(0034)/2017-EMR-I).

\medskip
\noindent
{\bf Conflict of Interests.} The authors declare that there is no conflict of interests 
regarding the publication of this paper.

\end{document}